\documentclass[11pt]{article}
\usepackage{amsfonts}
\usepackage{mathrsfs}
\usepackage{amsmath}
\usepackage{amssymb}
\usepackage{graphicx}
\usepackage{epic}
\usepackage{amsthm,latexsym}
\renewcommand{\paragraph}{\roman{paragraph}}
\def \T{\mathcal{T}}
\def \G{\mathcal{G}}
\def \U{\mathcal{U}}
\usepackage[]{caption2} 
\newtheorem{theorem}{\scshape \mdseries  Theorem}[section]
\newtheorem{lemma}[theorem]{\scshape \mdseries  Lemma}
\newtheorem{coro}[theorem]{\scshape \mdseries  Corollary}

\begin{document}

\title{\sf The spectral radius of the square of graphs\thanks{Supported by National Natural Science Foundation of China (11071002, 11371028),
Program for New Century Excellent Talents in University (NCET-10-0001),
Scientific Research Fund for Fostering Distinguished Young Scholars of Anhui University(KJJQ1001),
Academic Innovation Team of Anhui University Project (KJTD001B).}}
\author{Yi-Zheng Fan\thanks{Corresponding author. Email address: fanyz@ahu.edu.cn(Y.-Z. Fan), wanglongxuzhou@126.com(L. Wang)}, ~~Long Wang\\
  {\small  \it School of Mathematical Sciences, Anhui University, Hefei 230601, P. R. China}
 }
\date{}
\maketitle

\noindent {\bf Abstract:}  The square of a connected graph $G$ is obtained from $G$ by adding an edge between every pair of vertices at distance $2$.
In this paper we give some upper or lower bounds for the spectral radius of the square of connected graphs, trees and unicyclic graphs respectively.
We also investigate the spectral radius of the square of unicyclic graphs with given girth or trees with fixed diameter.

\noindent {\bf MR Subject Classifications:} 05C50


\section{Introduction}
Throughout this paper we consider simple, undirected and connected graphs.
Let $G=(V(G),E(G))$ be a connected graph of order $n$ and let $A(G)$ be its adjacent matrix.
Since $A(G)$ is a symmetric and nonnegative matrix, the largest eigenvalue of $A(G)$ is exactly the spectral radius of $A(G)$, which is also called the {\it spectral radius} of $G$ and is denoted by $\rho(G)$.

%Since the use of spectral radius $\rho(G)$ as a molecular descriptor was proposed by Balaban et al. (\cite{bal}), it has been well studied during recent years.  These are some basic famous results:
The following are among of the classical results of the spectral radius of graphs, where $P_n,C_n,S_n,K_n$ are respectively the path, the cycle, the star and the complete graphs of order $n$, and $S_{n}^{*}$ is obtained from $S_{n}$ by joining two of its pendant vertices.

\begin{theorem}\label{original} {\em \cite{bru,cve2, cve3, hon, li,simc}}
 Let $G$ be a connected graph of order $n$.
 If $G\neq P_{n}$ and $G\neq K_{n}$, then
$$\rho(P_{n})<\rho(G)<\rho(K_{n}).$$
If $G$ is a tree,  $G\neq P_{n}$ and $G\neq S_{n}$, then
$$\rho(P_{n})<\rho(G)<\rho(S_{n}).$$
If $G$ is a unicyclic graph, $G\neq C_{n}$ and $G\neq S_{n}^{*}$, then
$$\rho(C_{n})<\rho(G)<\rho(S_{n}^{*}).$$
\end{theorem}

The $k$-th {\it power} of a  connected graph $G$, denoted by $G^{k}$, is obtained from $G$ by adding an edge between every pair of vertices within distance $k$.
It is easy to see that $G^{1}=G$ and $G^{d}=K_{n}$ if $d$ is the diameter of $G$.
In particular, the graph $G^2$ is called the {\it square} of $G$.
If $G$ is a path or a cycle, then $G^{k}$ can also be considered as a distance graph or circulant graph.
   The power of graphs has been found to be useful in practical applications.
   For example, coloring on the square of a graph may be used to assign frequencies to the participants of wireless communication networks so that no two participants interfere with each other at any of their common neighbors \cite{agn,miao}, and to find graph drawings with high angular resolution \cite{for}.
   If $G$ is connected, then $G^{3}$ necessarily contains a Hamiltonian cycle; and $G^{2}$ is always Hamiltonian when $G$ is $2$-vertex connected (see \cite{die}).

In 1973 Cvetkovi\'{c} characterized the spectrum of the total graph of a regular graph \cite{cve1}.
Recall that the total graph of a graph $G$ is exactly $S(G)^{2}$, where $S(G)$ is the subdivision of $G$.
 Except for the above work, few work appears on the spectra of the power of graphs.
 In this paper, we give some upper or lower bounds for the spectral radius of the square of connected graphs, trees and unicyclic graphs respectively.
We also investigate the spectral radius of the square of unicyclic graphs with given girth or trees with fixed diameter.

\section{Preliminaries}
Denote by $\G_n,\T_n,\U_n$ the classes of connected graphs, trees and unicyclic graphs of order $n$ respectively.
A graph $G$ is called {\it maximizing} (respectively, {\it minimizing}) in a class of graphs if the spectral radius of the square of such graph attains the maximum (respectively, minimum) among all graphs in the class.

Let $G$ be a graph and let $u,v\in V(G)$.
 Denote by $N_{G}(u)$ the set of neighbors of $u$ in $G$, and by $d_{G}(u)=|N_{G}(u)|$ the degree of $u$.
  The distance of $u$ and $v$ is denoted by $dist_{G}(u,v)$, and the diameter of $G$ is denoted and defined as $diam(G)=\max_{u,v}dist(u,v)$.
  % A graph $G$ is called bipartite, if $V(G)$ is a disjoint union of $V_{1}$ and $V_{2}$, where $|V_{1}|=k$ and $|V_{2}|=l$, such that two end vertices of each edge comes from distinct $V_{i}$. If $uv\in E(G)$ for every $u\in V_{1}$ and $v\in V_{2}$, then $G$ is a complete bipartite graph and is denoted as $K_{k,l}$.   Assume $|V(G)|=n$ and $|E(G)|=m$, we call a connected graph $G$ $k$-cyclic if $m-n=k+1$ . If $k=0$, then $G$ is a tree.  One can see that $K_{n-1,1}$ is a tree of order $n$ with diameter $2$, it is called a star and also denoted as $S_{n}$. The vertex with degree $n-1$ is called the center of $S_{n}$. $S_{n}^{*}$ is obtained from $S_{n}$ by joining two pendant vertices. The double star $DS_{k,l}$ is a tree with diameter $3$, which is obtained from two stars $S_{k}$ and $S_{l}$ by identifying one pendant vertex.
The graph $G$ is called {\it regular} if all vertices of $G$ have the same degree, and is called {\it nontrivial} if $G$ contains more than one vertex.

A path on $n$ vertices is denoted by $P_n=v_1v_2\ldots v_n$, which contains the vertices $v_1,v_2,\ldots,v_n$ and edges $v_{i}v_{i+1}$ for $i=1,2,\ldots,n-1$.
A cycle on $n$ vertices is denoted by $C_n=v_1v_2\ldots v_n$, which contains the vertices $v_1,v_2,\ldots,v_n$ and edges $v_{i}v_{i+1}$ for $i=1,2,\ldots,n$ under modulo $n$.
The notation $H\subseteq G$ or $G \supseteq H$ means that $H$ is a subgraph of $G$, and
$H\subsetneq G$ or $G \supsetneq H$ means that $H$ is a proper subgraph of $G$ (i.e. $E(G) \backslash E(H) \ne \emptyset$).

\begin{lemma}\label{degree}{\em \cite{bro}} Let $G$ be a graph with maximum degree $\Delta(G)$
and average degree $\alpha(G)$.  If $G$ is not regular, then $\alpha(G)< \rho(G)< \Delta(G)$.
\end{lemma}

%In the If $H$ is a subgraph $H$ of $G$, we will denote is a graph satisfying $V(H)\subset V(G)$ and $E(H)\subset E(G)$, denoted as. If $\emptyset\neq E(H)\subsetneqq E(G)$, then $H$ is called a nontrivial subgraph of $G$, denoted as $H\subsetneqq G$. Furthermore, if $uv\in E(H)$ if and only if $uv\in E(G)$, for $u,v\in V(H)$, then $H$ is called an induced subgraph of $G$, denoted as $H\lhd G$.

\begin{lemma}\label{subgraph}{\em\cite{bro}} Let $G$ be a connected graph and $H$ be a subgraph of $G$.
 Then $\rho(H)\leq \rho(G)$, with equality if and only if $H=G$.
 \end{lemma}

Let $G_{1}$ and $G_{2}$ be two connected vertex-disjoint graphs, and $v_{1}\in G_{1}$, $v_{2}\in G_{2}$.
The {\it coalescence} of $G_{1}$ and $G_{2}$, denoted by $G_{1}(v_{1})\circ G_{2}(v_{2})$, is obtained from $G_{1}$ and $G_{2}$ by identifying $v_{1}$ with $v_{2}$ and forming a new vertex $u$.
The graph $G_{1}(v_{1})\circ G_{2}(v_{2})$ is also written as $G_{1}(u)\circ G_{2}(u)$.  If a connected graph $G$ can be expressed in the form $G=G_{1}(u)\circ G_{2}(u)$, where $G_{1}$ and $G_{2}$ are both nontrivial and connected, then $G_{1}$ and $G_{2}$ are called the {\it branches} of $G$ rooted at $u$.
If $v'_1$ is another vertex of $G_1$ rather than $v_1$, the graph $G_{1}(v'_{1})\circ G_{2}(v_{2})$ is said obtained from $G_{1}(v_{1})\circ G_{2}(v_{2})$ by {\it relocating} $G_2$ from $v_1$ to $v'_1$.

\begin{lemma}\label{coale}  Let $G$ be a nontrivial connected graph and $T$ be a nontrivial tree.
 Then $$\rho[(G(v)\circ T(u))^{2}]\leq \rho[(G(v)\circ S(u))^{2}],$$ where $S$ is a star centered at $u$  that has the same order as $T$.
   The equality holds if and only if $T$ is a star centered at $u$.
   \end{lemma}

\noindent{\bf Proof.} It can be directly checked that, if $T$ is not a star centered at $u$,
then $G(v)\circ T(u)\subsetneq G(v)\circ S(u)$. The result follows from Lemma \ref{subgraph}.
\hfill   $\blacksquare$

\begin{coro} \label{quasi}
Let $P_3=v_{1}v_{2}v_{3}$ and let $H$ be a nontrivial connected graph.
Let $G_{1}=P_{3}(v_{3})\circ H(u)$ and $G_{2}=P_{3}(v_{2})\circ H(u)$.
Then $\rho(G_{1}^{2})<\rho(G_{2}^{2})$.
\end{coro}

From the famous Perron-Frobenius theory, for a connected graph $G$ of order $n$, $\rho(G^{2})$ is simple and is corresponding to a positive eigenvector, which is called the {\it Perron vector} of $G^{2}$. Let $X$ be a Perron vector of $G^{2}$, and let $X_{v}$ denote the entry of $X$ corresponding to the vertex $v$ of $G$. We have
$$\rho(G^{2})X_{v}=\sum_{u,dist_{G}(u,v)\leq 2}X_{u}.\eqno(2.1)$$
Denote $\tilde{X_{v}}=X_{v}+\sum_{u\in N_{G}(v)}X_{u}$, we have the following useful lemma.

\begin{lemma}\label{pendant} Let $H_{1}$ and $H_{2}$ be two nontrivial connected graphs,  and let $u,v\in V(H_{1})$, $w\in V_{H_{2}}$.
Let $G_{1}=H_{1}(u)\circ H_{2}(w)$ and $G_{2}=H_{1}(v)\circ H_{2}(w)$.
Let $X$ be a Perron vector of $G_{1}^{2}$.
If $\tilde{X_{u}}\leq \tilde{X_{v}}$ and $X_{u}\leq X_{v}$, then $\rho(G_{1}^{2})<\rho(G_{2}^{2})$.
In particular, if  $H_{2}$ is a star, and $\tilde{X_{u}}\leq \tilde{X_{v}}$, then $\rho(G_{1}^{2})<\rho(G_{2}^{2})$.\end{lemma}

\noindent{\bf Proof.}  Suppose that $X$ is unit.
We have
\begin{align*}
    \frac{1}{2}X^{T}A(G_{2}^{2})X &=\sum_{rs\in E(G_{2}^{2})}X_{r}X_{s} \\
    &=\sum_{rs\in E(G_{1}^{2})}X_{r}X_{s}+\sum_{r\in N_{H_{2}(w)}}\sum_{s\in N_{H_{1}}(v)\cup \{v\}}X_{r}X_{s}
    +\sum_{r,dist_{H_{2}}(r,w)=2}X_{r}X_{v}\\
    &~~-\sum_{r\in N_{H_{2}}(w)}\sum_{s\in N_{H_{1}}(u)\cup \{u\}}X_{r}X_{s}-\sum_{r,dist_{H_{2}}(r,w)=2}X_{r}X_{u}\\
    &=\frac{1}{2}X^{T}A(G_{1}^{2})X+\sum_{r\in N_{H_{2}}(w)}X_{r}(\tilde{X_{v}}-\tilde{X_{u}})
      +\sum_{r,dist_{H_{2}}(r,w)=2}X_{r}(X_{v}-X_{u})\\
    &\geq \frac{1}{2}X^{T}A(G_{1}^{2})X.
\end{align*}
So $$\rho(G_{2}^{2})\ge X^{T}A(G_{2}^{2})X \ge X^{T}A(G_{1}^{2})X =\rho(G_{1}^{2}).$$

If $\rho(G_{1}^{2})=\rho(G_{2}^{2})$, then $X$ is also the Perron vector of $G_{2}^{2}$.
Furthermore, $\tilde{X_{v}}=\tilde{X_{u}}$, and $X_{v}=X_{u}$ if there exists a vertex $r$ such that $d_{H_{2}}(r,w)=2$.
However, by considering the eigenvector equation (2.1) for $G_1^2$ and $G_2^2$ both on the vertex $u$,
we get that $u$ is adjacent to $v$, and $H_2$ is a star centered at $w$.
As $\tilde{X_{v}}=\tilde{X_{u}}$, there must exist a vertex $t$ of $H_1$ other than $u$ and $v$, such that $t$ is adjacent to $v$ but not adjacent to $u$.
Note that $N_{G_2}(t) \backslash N_{G_1}(t)=V(H_2)\backslash \{w\}$.
So, if considering the  equations (2.1) for $G_1^2$ and $G_2^2$ both on the vertex $t$, we also get a contradiction.

The second result is easily obtained by the above discussion. \hfill   $\blacksquare$

\begin{coro} \label{rotate}
Let $T$ be a tree obtained from $P_4=v_1v_2v_3v_4$ by appending a pendant edge $uv_3$ at $v_3$, let $H$ be a connected graph containing a vertex $w$.
Let $G_1=T(v_4) \circ H(w)$ and let $G_2$ be obtained from $G_1$ by relocating the edge $uv_3$ to $uv_1$.
Then $\rho(G_2^2) < \rho(G_1^2)$.
\end{coro}

\noindent{\bf Proof.}
Let $X$ be a Perron vector $G_2^2$.
If $\tilde{X}_{v_1}\leq \tilde{X}_{v_3}$, then $\rho(G_2^2)<\rho(G_1^2)$ by Lemma \ref{pendant}.
Otherwise, $\tilde{X}_{v_1}> \tilde{X}_{v_3}$, i.e. $X_{u}+X_{v_{1}}+X_{v_{2}} > X_{v_{2}}+X_{v_{3}}+X_{v_{4}}$.
If $X_{v_1} < X_{v_3}$, then $X_u > X_{v_4}$.
However, by the assumption $X_{v_1} < X_{v_3}$ and the eigenvector equation (2.1) on $u$ and $v_4$, we get $X_u < X_{v_4}$, a contradiction.
So, $X_{v_1} \ge X_{v_3}$.

Observe that $G_2=P(v_3) \circ (H+v_3v_4)(v_3)$, where $P=uv_1v_2v_3$ and $H+v_3v_4$ is obtained from $H$ by appending a pendant edge $v_3v_4$.
Now relocating the branch $G+v_3v_4$ from $v_3$ to $v_1$, the resulting graph $\bar{G}_1$ is isomorphic to $G_1$.
By Lemma \ref{pendant}, we have $\rho(G_2^2)<\rho(\bar{G}_1^2)=\rho(G_1^2)$.
\hfill   $\blacksquare$

\section{The spectral radius of the square of graphs or trees}
In this section we will determine the maximizing or minimizing graphs in $\G_n$ or $\T_n$.
It is easily seen that $G^2$ is a complete graph if and only if $diam(G) \le 2$.
So, by Lemma \ref{subgraph} we have the following result.

\begin{theorem}\label{graph-large}
Let $G$ be a connected graph of order $n$.
Then $\rho(G^{2})\leq n-1$, with equality holds if and only if $diam(G) \le 2$, or equivalently $G^2=K_n$.
\end{theorem}

%From Theorem \ref{graph-large}, we know that the star $S_n$ is the unique tree whose square has the maximum spectral radius among all trees of order $n$.
We now wish to find the minimizing graph in $\G_n$. By Lemma \ref{subgraph}, such graph can be found within $\T_n$.
Let $T_{i}$ be the graphs listed in Fig. \ref{pic-tree}. By {\scshape \mdseries  Mathematica}, $\rho(T_{i}^{2})> 4$ for $i=1,2, \ldots,8$.
Note that $\rho(P_{n}^{2})< \Delta(P_{n}^{2})=4$.
So, if $T$ is a minimizing graph in $\T_n$, then $T$ cannot contain any graph in Fig. \ref{pic-tree} as a subgraph.
In other words, any graph in Fig. \ref{pic-tree} is {\it forbidden} in $T$.
%In particular, from the graph $T_1$, we find that $\Delta(T) \le 3$.

\renewcommand\thefigure{\arabic{section}.\arabic{figure}}
\begin{figure}[h!]
  \centering
  \includegraphics[scale=0.6]{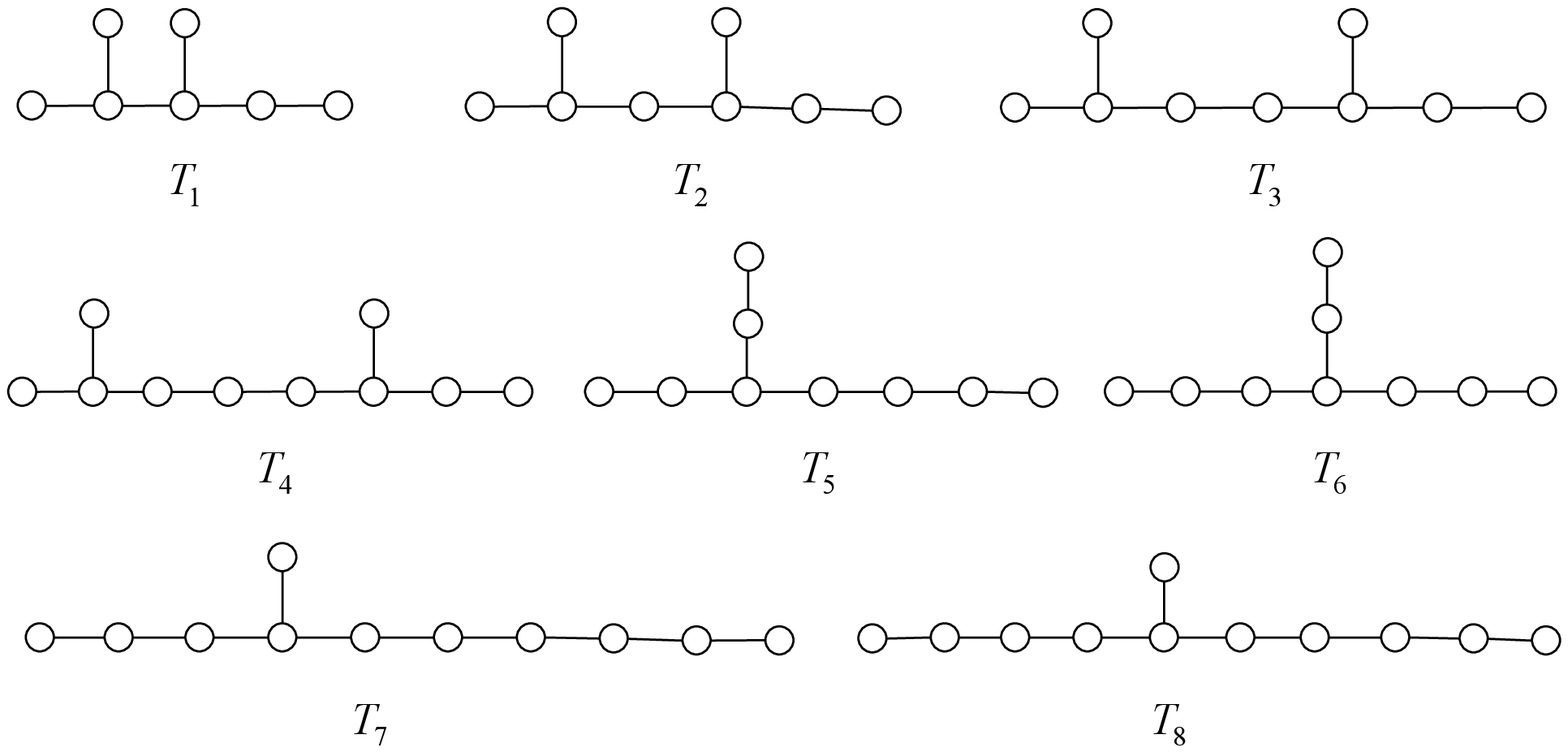}\\
  \caption{\small Eight trees $T_i$ with $\rho(T_i^2)>4$ for $i=1,2,\ldots,8$}\label{pic-tree}
\end{figure}

\begin{lemma}\label{paths}
Let $P_{k}=v_{1}v_{2}\ldots v_{k}$ and $P_{n-k+1}=u_{1}u_{2}\ldots u_{n-k+1}$ be two disjoint paths.
Assume that $n\geq 2k-1$ and $T=P_{k}(v_{k})\circ P_{n-k+1}(u_{k})$.
Then $\rho(T^{2})\geq \rho(P_{n}^{2})$, with equality if and only if $n=2k-1$.
\end{lemma}

\noindent{\bf Proof.}
If $n=2k-1$, clearly the result follows as $T$ is exactly a path of order $n$.
Now suppose that $n >2k-1$.
 Let $P_{n}=w_{1}w_{2}\ldots w_{n}$, $P_{n}^{(1)}=w_{1}w_{2}\ldots w_{n-k+1}$ and $P_{n}^{(2)}=w_{n-k+1}w_{n-k+2}\ldots w_n$ be its two subpaths.
 One can observe that $\bar{T}:=P_{n}^{(2)}(w_{n-k+1}) \circ P_{n}^{(1)}(w_k)$ is isomorphic to $T$.
 Let $X$ be the Perron vector of $P_{n}^{2}$.
 Then by symmetry $\tilde{X}_{w_{k}}=\tilde{X}_{w_{n-k+1}}$ and $X_{w_{k}}=X_{w_{n-k+1}}$.
 So, by Lemma \ref{pendant}, $\rho(P_{n}^{2})< \rho(\bar{T}^{2})=\rho(T^{2})$.\hfill   $\blacksquare$

\begin{lemma}\label{v3} Let $T$ be a tree of order $n$ with diameter at least $4$, and
 let $P_{k}=v_{1}v_{2}\ldots v_{k}$ be a longest path in $T$.
  If $T$ is a minimizing graph in $\T_n$, then  $d_{T}(v_{2})=d_{T}(v_{3})=d_{T}(v_{k-1})=d_{T}(v_{k-2})=2$.
  \end{lemma}

\noindent{\bf Proof.} By Corollary \ref{quasi} we find that $d_{T}(v_{2})=d_{T}(v_{k-1})=2$.
First we assert that, if $d_{T}(v_{3})\geq 3$, then $d_T(u)=1$ for any $u\in N_{T}(v_{3})\backslash V(P_{k})$. (In fact, if $N_{T}(v_{3})\backslash V(P_{k}) \ne \emptyset$,
then it contains only one element as $S_5$ is forbidden in $T$.)
 Otherwise, if $k\geq 7$, then $T \supseteq T_{5}$, and hence $\rho(T^2) \ge \rho(T_5^{2})>4$, a contradiction, where $T_{5}$ is shown in Fig. \ref{pic-tree}.
 So, $k$ takes $5$ or $6$.
 As $T$ contains no $T_{1}$ in Fig. \ref{pic-tree} as a subgraph, $T=P_{3}(v_{3})\circ P_{n-2}(u_{3})$, where $P_3$ and $P_{n-2}$ are denoted as in Lemma \ref{paths},
 and $n$ takes  $7$ or $8$.
 By Lemma \ref{paths}, we get $\rho(T^{2})>\rho(P_{n}^{2})$, a contradiction.

Now suppose that $N_{T}(v_{3})\backslash V(P_{k})$ contains only one pendant vertex, say $u$.
%Observe that $T=T_{(1)}(v_{3})\circ T_{(2)}(v_{3})$, where $T_{(1)}$ is a tree induced on the vertices $u$ and $v_3$.
Let $\bar{T}$ be obtained from $T$ by relocating the edge $uv_3$ to $uv_1$.
%, and let $X$ be a Perron vector $\bar{T}^2$.
%If $\tilde{X}_{v_1}\leq \tilde{X}_{v_3}$, then $\rho(\bar{T}^2)<\rho(T^{2})$ by Lemma \ref{pendant}, a contradiction.
%Otherwise, $\tilde{X}_{v_1}> \tilde{X}_{v_3}$, i.e. $X_{u}+X_{v_{1}}+X_{v_{2}} > X_{v_{2}}+X_{v_{3}}+X_{v_{4}}$.
%If $X_{v_1} < X_{v_3}$, then $X_u > X_{v_4}$.
%However, by the assumption $X_{v_1} < X_{v_3}$ and the eigenvector equation (2.1) on $u$ and $v_4$, we get $X_u < X_{v_4}$, a contradiction.
%So, $X_{v_1} \ge X_{v_3}$.
%We also observe that $\bar{T}=\bar{T}_1(v_3) \circ \bar{T}_2(v_3)$, where $\bar{T}_1$ is induced on the vertices $u,v_1,v_2,v_3$.
%Now relocating the branch $\bar{T}_2$ from $v_3$ to $v_1$, the resulting graph $T'$ is isomorphic to $T$.
By Corollary \ref{rotate} we get $\rho(\bar{T}^2)< \rho(T^2)$, a contradiction.

By the above discussion, we have shown $d_{T}(v_{3})=2$. By symmetry we also have $d_{T}(v_{k-2})=2$.\hfill   $\blacksquare$

\begin{theorem}\label{tree-small}
Let $T$ be a tree of order $n\geq 4$.
Then $\rho(T^{2})\geq \rho(P_{n}^{2})$, with equality holds if and only if $T=P_{n}$.
\end{theorem}

\noindent{\bf Proof.} Let $T$ be a minimizing tree in $\T_n$.
We assert that $T=P_n$ and hence the result follows.
Let $P_{k}=v_{1}v_{2}\ldots v_{k}$ be a longest path in $T$.
The result holds if $k$ equals $3$ or $4$ by Corollary \ref{quasi}.
If $k \ge 5$, then $d_T(v_2)=d_{T}(v_{3})=d_{T}(v_{k-2})=d_{T}(v_{k-1})=2$ by Lemma \ref{v3}.
So the result holds for $k$ being $5$ or $6$.

Now suppose $k \ge 7$.
Assume to the contrary that $T \ne P_n$.
If $n \geq k+2$, then $T$ contains one of $T_{1},T_{2},T_{6},T_{7}$ in Fig. \ref{pic-tree} as a subgraph, which implies that $\rho(T^2)>4$, a contradiction.

If $n=k+1$, letting $u\in V(T)\backslash V(P_{k})$, we assert that if $uv_{4}\in E(T)$ or $uv_{v_{k-3}}\in E(T)$, then $\rho(T^{2})>\rho(P_{n}^{2})$.
Assume without loss of generality that $uv_{4}\in E(T)$.
Let $P_{n}=u_{1}u_{2}\ldots u_{n}$ be the path of order $n$ and let $X$ be the Perron vector of $P_{n}^{2}$.
Then $P_n=P_n^{(11)}(v_{n-1}) \circ P_n^{(12)}(v_{n-1})$, where $P_n^{(12)}=u_{n-1}u_n$; or
$P_n=P_n^{(21)}(u_{4}) \circ P_n^{(22)}(u_{4})$, where $P_n^{(21)}=u_1u_2u_3u_4$.
If $\tilde{X}_{u_{4}} \geq \tilde{X}_{u_{n-1}}$, noting that $T$ is obtained from $P_n$ by relocating $P_n^{(12)}$ from $u_{n-1}$ to $u_4$, hence by Lemma \ref{pendant} we have $\rho(P_{n}^{2})<\rho(T^{2})$.
Otherwise, $\tilde{X}_{u_{4}} < \tilde{X}_{u_{n-1}}=\tilde{X}_{u_2}$, and then $X_{u_4}+X_{u_5}<X_{u_1}+X_{u_2}$.
If $X_{u_4}>X_{u_2}$, then $X_{u_5}<X_{u_1}$.
However, by the assumption that $X_{u_4}>X_{u_2}$ and the eigenvector equation (2.1) on $u_1$ and $u_5$, we get $X_{u_5}>X_{u_1}$, a contradiction.
So $X_{u_4}\le X_{u_2}=X_{u_{n-1}}$.
Noting that $T$ is also obtained from $P_{n}$ by relocating $P_n^{(21)}$ from $u_4$ to $u_{n-1}$,
so by Lemma \ref{pendant} again we get $\rho(P_{n}^{2})<\rho(T^{2})$.

From the above discussion, it suffices to consider the case of $k \ge 9$.
If $k=9$, then $\rho(T^{2})>\rho(P_{10}^{2})$ by Lemma \ref{paths}.
If $k \ge 10$, then $T\supseteq T_{8}$ and $\rho(T^{2})>4$.
The result follows now.\hfill   $\blacksquare$

\begin{coro}\label{tree-large}
Let $T$ be a tree of order $n\ge 4$.
Then $$\rho(P_{n}^{2})\leq \rho(T^{2})\leq \rho(S_{n}^{2}),$$
 with the first equality if and only if $T=P_{n}$, and the second equality if and only if $T=S_{n}$.
 \end{coro}

\begin{coro}
Let $G$ be a connected graph of order $n\geq 4$.
Then $\rho(G^{2})\geq \rho(P_{n}^{2})$, with equality holds if and only if $G=P_{n}$.
\end{coro}

\noindent{\bf Proof.}
If $G$ itself is a tree, the result follows from Theorem \ref{tree-small}.
Otherwise, $G$ contains a spanning unicyclic subgraph $U$ which contains a cycle $C$.
Surely by Lemma \ref{subgraph}, $\rho(G^2) \ge \rho(U^2)$.
If $U=C$, noting $U$ contains at least $4$ vertices, we have $\rho((U-e)^2)<\rho(U^2)$ as $(U-e)^2$ is a proper subgraph of $U^2$, where $e$ is an arbitrary edge on the cycle.
Otherwise, let $uv$ be an edge appending on the cycle $C$ and let $uw=:e$ be an edge on $C$, where $u \in V(C)$.
Now the distance between $v,w$ is $2$ in $U$, but is greater than $2$ in $U-e$.
So, $(U-e)^2$ is still a proper subgraph of $U^2$, and hence $\rho((U-e)^2)<\rho(U^2)$.
By the above discussion, $$\rho(G^2) \ge \rho(U^2)>\rho((U-e)^2) \ge \rho(P_{n}^{2}).$$
\hfill   $\blacksquare$

%Until now we have concluded that $P_{n}$ is the unique graph having the smallest square spectral radius among trees. To complete the proof, we only need to show that $G^{2}=P_{n}^{2}$ if and only if $G=P_{n}$ when $n\geq 4$. In fact, If $G$ contains at least one cycle, then $G\supset C_{k}$ for some $k\geq 5$, which leads $\rho(G^{2})\geq \rho(C_{k}^{2})=4$; or $G$ contains $S_{4}^{*}$, which leads $G^{2}\supset K_{4}$, while $K_{4}\not\subset P_{n}^{2}$. \hfill   $ \Box $ \\

\noindent{\bf Remark.} The result dose not hold for $n=3$ as $\rho(P_{3}^{2})=\rho(C_{3}^{2})=2$.

%Note that  for a tree $T$, $T^{2}=K_{n}$ if and only if $T=S_{n}$, and the tree with the smallest  square spectral radius has been determined in Theorem \ref{graph-small}, we have the following theorem.

\section{The spectral radius of the square of unicyclic graphs}
In this section we determine the maximizing or minimizing graphs in $\U_n$.
The following result on maximizing graphs in $\U_n$ is easily seen, that is, the graphs $C_{3},C_{4},C_{5},S_{n}^{*}$ with diameter at most $2$.

\begin{theorem}\label{uni-large} Let $U$ be a unicyclic graph of order $n\geq 3$.
 Then $\rho(U^{2})\leq n-1$, with equality if and only if $diam(U) \le 2$.
 \end{theorem}

Now we characterize the minimizing graph in $\U_n$.

\begin{lemma}\label{uni-small-girth} Let $U$ be a unicyclic graph of order $n \ge 6$.
If $5\leq g(U)\leq n-1$, then $\rho(U^{2})> 4$.
\end{lemma}

\noindent{\bf Proof.}
By Lemma \ref{degree}, it suffices to prove that the average degree $\alpha(U^2)$ of $U^2$ is greater than $4$.
Denote $g:=g(U)$.
If $n=g+1$, then $\alpha(U^2) >4$ obviously by a direct calculation.
Assume that $\alpha(U^2) >4$ for all unicyclic graphs $U$ of order $n = g+k$, where $k \ge 1$.
Let $U$ be unicyclic graph of order $g+k+1$.
Then $U$ is obtained from a unicyclic graph, say $\bar{U}$ of order $g+k$ by appending a pendant edge $uw$, where $w$ is the vertex of degree one.
Note that $u$ has at least one neighbor $u'$ in $\bar{U}$.
So, $d_{U^2}(w)\ge 2$, $d_{U^2}(u)= d_{\bar{U}^2}(u)+1$ and $d_{U^2}(u')= d_{\bar{U}^2}(u')+1$.
By the assumption,
\begin{align*}
\sum_{v \in U}d_{U^2}(v) & =  \sum_{v \in U\backslash \{u,u'\}}d_{U^2}(v)+d_{U^2}(u)+d_{U^2}(u')+d_{U^2}(w) \\
&\ge \sum_{v \in U\backslash \{u,u'\}}d_{U^2}(v)+(d_{\bar{U}^2}(u)+1)+(d_{\bar{U}^2}(u')+1)+2\\
& \ge \sum_{v \in U\backslash \{u,u'\}}d_{\bar{U}^2}(v)+d_{\bar{U}^2}(u)+d_{\bar{U}^2}(u')+4\\
&= \sum_{v \in \bar{U}}d_{\bar{U}^2}(v)+4\\
& = (g+k)\alpha(\bar{U})+4\\
& > 4(g+k+1)
\end{align*}\hfill   $\blacksquare$
%If $g(U)=k\geq 5$, we prove $\alpha(U^{2})\geq 4$ by induction on $n$ and the result will follow from Lemma \ref{degree}.  If $n=k$, then $\alpha(U^{2})=4$ trivially. Assume by induction that $\alpha(U^{2})\geq 4$ when $n\leq k+l$.
%
%Now consider an arbitrary unicyclic graph $\tilde{U}$ with girth $k$ and order $n=k+l+1$. $\tilde{U}$ can be viewed as obtained from some unicyclic graph $\bar{U}$ (with girth $k$ and order $k+l$) by adding a pendant vertex $v$. Denote $N_{\tilde{U}^{2}}(v)=\{u_{1},u_{2},\ldots,u_{k}\}$, then $k\geq 2$. Note that $d_{\tilde{U}^{2}}(u_{i})\geq d_{\bar{U}^{2}}(u_{i})+1$ for $1\leq i\leq k$, we have $\sum_{j=1}^{k+l+1}d_{\tilde{U}^{2}}(v_{j})\geq \sum_{j=1}^{k+l}d_{\bar{U}^{2}}(v_{j})+4$ and $\alpha(\tilde{U}^{2})\geq\alpha(\bar{U^{2}})\geq4$. \hfill   $ \Box $ \\

Note that $\rho(C_{n}^{2})=4$ when $n\geq 5$.
Let $U$ be  a unicycle graph of order $n$.
The above lemma shows that if $U$ is not a cycle and $g(U) \ge 5$, then  $U$ is not a minimizing graph.
Next we consider the unicycle graphs with girth less than $5$.
By a calculation using {\mdseries \sc Mathematica}, the graphs in Fig. \ref{pic-uni-small} are forbidden in a minimizing graph.

% the above lemma shows if $g(U)\geq 5$ and $U\neq C_{n}$, then $\rho(U^{2})$ can not attain the smallest square spectral radius among unicyclic graphs.
\setcounter{figure}{0}
\renewcommand\thefigure{\arabic{section}.\arabic{figure}}
\begin{figure}[h!]
  \centering
  \includegraphics[scale=0.6]{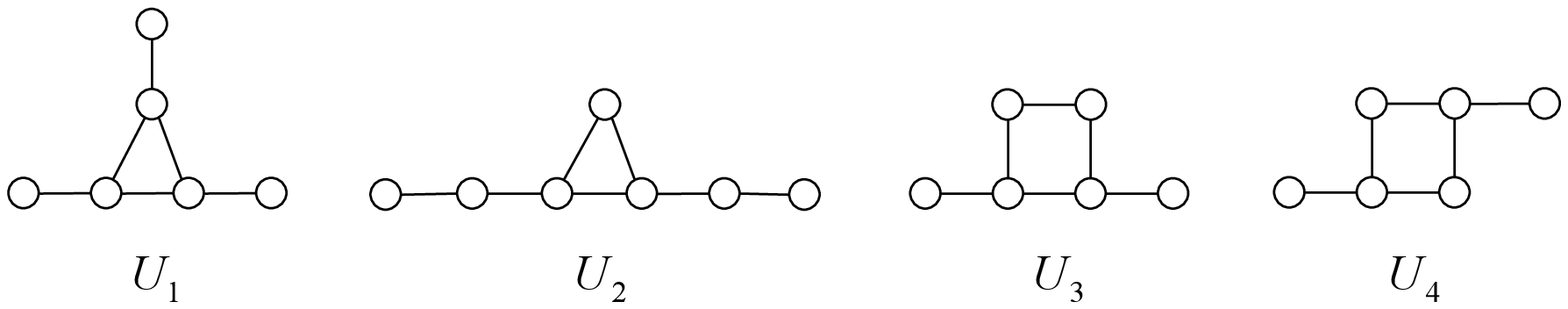}\\
  \caption{\small Four unicyclic graphs $U_i$ with $\rho(U_i^2)>4$ for $i=1,2,3,4$}\label{pic-uni-small}
\end{figure}

\begin{lemma}\label{uni-small-tri}
Let $U$ be a minimizing graph in $\U_n\;(n \ge 5)$.
If $g(U)=3$ or $g(U)=4$, then $U=C_{3}(v)\circ T(v)$, where $T$ is a tree containing $v$ as a pendant vertex.
\end{lemma}

\noindent{\bf Proof.}
We first consider the case of $g(U)=3$.
Let $C_3=v_1v_2v_3$ be the cycle of $U$.
As $U_{1}$ of Fig. \ref{pic-uni-small} is forbidden, $C_3$ has at most two vertices of degree at least $3$;
If $C_3$ contains exactly two vertices, say $v_1,v_2$ of degree at least $3$, then one of them, say $v_1$ is only adjacent to pendant vertices as $U_2$ of  Fig. \ref{pic-uni-small} is forbidden.
Furthermore, by Corollay \ref{quasi}, the degree of $v_1$ is exactly $3$.
Let $u$ be the unique neighbor of $v_1$ outside $C_3$.
Now delete the edge $v_{2}v_{3}$ and adding a new edge $uv_{3}$, we will get a graph $U'$ which holds $U'^{2}\subsetneq U^{2}$.
So $\rho(U'^{2})<\rho(U^{2})$ by Lemma \ref{subgraph}, a contradiction to $U$ be minimizing.

By the above discussion, $C_3$ contains exactly one vertex, namely $v_2$ of degree at least $3$.
If $d_U(v_2)=3$, we are done.
Otherwise, $U$ contains $S_5^*$ as a proper subgraph or $U=S_5^*$.
The latter case cannot happen by Corollay \ref{quasi}.
For the former case, $\rho(U^2)>4$, also a contradiction.

Now come to the case of $g(U)=4$.
Let $C_4=v_1v_2c_3v_4$ be the cycle of $U$.
First suppose that $n \ge 6$.
As $U_{3}$ and $U_{4}$ of Fig. \ref{pic-uni-small} are forbidden, $C_4$ contains exactly one vertex, say $v_4$ of degree at least $3$.
In fact, $d_U(v_4)=3$; otherwise $U$ would contain a forbidden subgraph $S_5$.
So we can write $U$ as $U=C_{4}(v_4)\circ T(v_4)$, where $T$ is a tree containing $v_4$ as a pendant vertex.
If $n=5$, then $U$ is obtained from $C_4$ by a appending a pendant edge at some vertex say $v_4$.
Now construct a new unicyclic graph $\bar{U}$ from $U$ by deleting the edge $v_{3}v_4$ and adding a new edge $v_{3}v_{1}$.
Then $U^{2}\supsetneq \bar{U}^{2}$ and $\rho(U^{2})>\rho(\bar{U}^{2})$, a contradiction. \hfill   $\blacksquare$

\begin{theorem}\label{uni-small}
Let $U$ be a unicyclic graph of order $n \ge 4$.
Then $$\rho(U^{2})\geq \min\{\rho[(C_{3})(v)\circ P_{n-2}(v))^{2}],\rho(C_{n}^{2})\},$$
 with equality if and only if $U=C_{3}(v)\circ P_{n-2}(v)$ or $U=C_{n}$, where $v$ is a pendant vertex of $P_{n-2}$.
 \end{theorem}

\noindent{\bf Proof.}  Surely the result holds for $n=4$.
 Suppose $n \ge 5$ in the following.
 Assume that $U \ne C_n$ is a minimizing graph in $\U_n$.
By Lemma \ref{uni-small-girth} and Lemma \ref{uni-small-tri}, $U=C_{3}(v_1)\circ T(v_1)$, where $T$ is a tree containing $v_1$ as a pendant vertex, and $d_U(v_1)=3$.
Let $P_{k}=v_{1}v_{2}\ldots v_{k}$ be the longest path contained in $T$ which starts from $v_1$, where $k \ge 3$ as $n \ge 5$.
 We assert that $U=C_{3}(v_1)\circ P_{n-2}(v_1)$ by the following claims.

(1) $d_U(v) \le 3$ for any vertex $v \in V(T)$; otherwise $U$ would properly contain $S_5$.

(2) $d_U(v_{k-1}) =2$; otherwise, the neighbors of $v_{k-1}$ outside $P_k$ should be pendant vertices by the definition of $P_k$, and
then $U$ is not minimizing by Corollary \ref{quasi}.

(3) $d_U(v_{k-2})=2$ by Corollay \ref{rotate}.

(4) $d_U(v_2)=d_U(v_3)=2$; otherwise, by the fact that $d_U(v_{k-1}) =2$ we find that $U$ contains $T_1$ or $T_2$ of Fig. \ref{pic-tree}.

(5) $d_U(v_4)=d_U(v_5)=2$; otherwise, $U$ contains $T_3$ or $T_4$ of Fig. \ref{pic-tree}.

%(5) for $i \ge 6$, the neighbors of $v_i$ should be  pendant vertices; otherwise, $U$ would contain $T_5$ of Fig. \ref{pic-tree}.

(6) for $6 \le i \le k-3$ (in this case $k \ge 9$), $d_U(v_i)=2$; otherwise $U$ would contain $T_7$  of Fig. \ref{pic-tree}.
\hfill   $\blacksquare$

%
%
%If some $u\notin \{v_{0},v_{1},\ldots,v_{k}\}$ is not pendant in $U$, then $U$ contains one of $\{T_{1},T_{2},T_{3},T_{5}\}$ as a subgraph and $\rho(U^{2})>4$. (Recall $T_{i}$ for $1\leq i\leq 13$ are shown in Fig. 1.) If $d_{U}(v_{k-1})\geq 3$, then $\rho(U^{2})$ is not the smallest follows from Lemma \ref{coale}. Assume $d_{U}(v_{k-1})$ but $u$ is pendant at $v_{k-2}$. Let $U'$ be graph obtained from $U$, by deleting $uv_{k-2}$ and adding $uv_{k}$. Then $U'$ contains a largest paths that starts at $C_{3}$, labelled as $u_{1}u_{2}\ldots u_{k}u_{k+1}$, where $u_{i}=v_{i}$ for $1\leq i\leq k$ and $u_{k+1}=u$. Note that $U$ can be obtained from $U'$ by deleting $u_{k+1}u_{k}$ and adding $u_{k+1}u_{k-2}$, also be obtained from $U'$ by deleting $u_{k-2}u_{k-3}$ and adding $u_{k}u_{k-3}$. Let $X$ be the Perron vector of $U'^{2}$, if $\tilde{X_{u_{k-2}}}\geq \tilde{X_{u_{k}}}$, then $\rho(U^{2})>\rho(U'^{2})$ by Lemma \ref{pendant}; If $\tilde{X_{u_{k-2}}}\leq \tilde{X_{u_{k}}}$, similar as in the proof of Lemma \ref{v3}, we have $X_{u_{2}}\leq X_{u_{k}}$, then $\rho(U^{2})>\rho(U'^{2})$ again.
%
%Now assume $u$ is pendant to $v_{i}$, where $3\leq i\leq k-3$. If $i\geq 7$, then $U\supset T_{6}$; If $3\leq i\leq 6$, then $U$ contains one of $\{T_{2},T_{3},T_{8},T_{9}\}$ as a subgraph ($T_{8}$ and $T_{9}$ are shown in Fig. 2, then $\rho(U^{2})>4$. \hfill   $ \Box $ \\

{\bf Conjecture 1.} $$\rho[(C_{3}(v)\circ P_{n-2}(v))^{2}]<\rho(C_{n}^{2}).$$
Using {\mdseries \sc Mathematica}, we verified the conjecture is true for $5\leq n\leq 100$.
If the conjecture were true, we would find some difference between the minimum spectral radius of unicyclic graphs (see Theorem \ref{original}) and their square.

\section{The spectral radius of the square of unicyclic graph with fixed girth}

Denote by $\U_n(g)$ the class of unicyclic graphs of order $n$ with girth $g$, and by
$\T_n(d)$ the class of trees of order $n$ with diameter $d$.
We determines the unique maximizing graph  in $\U_n(g)$, and then pose a conjecture on the maximizing graph(s) in $\T_n(d)$.

\begin{theorem}\label{uni-large-girth}
Let $\tilde{U}=C_{g}(v)\circ S_{n-g+1}(v)$, where $v$ is the center of $S_{n-k+1}$.
 Then for any unicyclic graph $U \in \U_n(g)$ of order $n$ with girth $g$, $\rho(U^{2})\leq \rho(\tilde{U}^{2})$, with the equality if and only if $U=\tilde{U}$.
 \end{theorem}

\noindent{\bf Proof.} Let $C_{g}=v_{1}v_{2}\ldots v_{g}$.
 Then $U$ can be considered as one obtained from $C_g$ by attaching some trees $T_{i}$ at the vertices $v_{i}$'s respectively.
 By Lemma \ref{coale}, we have $\rho(U^{2})\leq \rho(\bar{U}^{2})$, where $\bar{U}$ is obtained from $U$ by replacing each $T_i$ by a star centered $v_i$ that has the same order of $T_i$. Furthermore, the equality holds if and only if each $T_i$ is a star centered at $v_i$.

Suppose that there exists at least two stars of $\bar{U}$, say $S_{(i)}$ and $S_{(j)}$ centered at $v_i$ and $v_j$ respectively.
Let $X$ be a Perron vector of $\bar{U}^{2}$.
Assume without loss of generality that $\tilde{X}_{v_{i}}\leq \tilde{X}_{v_{j}}$.
 then $\rho(\bar{U}^{2})<\rho(\bar{U}'^{2})$ by Lemma \ref{pendant}, where $\bar{U}'$ is obtained from $\bar{U}$ by relocating $S_{(i)}$ from $v_i$  to $v_j$.
 Repeating the above process until the resulting graph contains only one vertex with degree greater than $2$, we then get the
 result as desired. \hfill   $\blacksquare$

%Now we prove $\rho(\bar{U}^{2})$ cannot attain the largest square spectral radius among unicyclic graphs with order $n$ and girth $k$, unless $\bar{U}=\tilde{U}$. Otherwise, assume $\bar{U}\neq \tilde{U}$ and $\rho(\bar{U}^{2})$ is the largest. Then at least two vertices $v_{i}$ and $v_{j}$ in $\bar{U}$ has degree at least $3$. Say $u_{1},u_{2},\ldots,u_{k}$ are pendant at $v_{i}$. Let $X$ be the

By a similar discussion as in the proof of Theorem \ref{uni-large-girth}, we narrow the scope the maximizing graph(s) in $\T_n(d)$.

\begin{theorem}
Let $P_{k+1}=v_{1}v_{2}\ldots v_{d+1}$, and let $\tilde{T}(i)=P_{k+1}(v_i)\circ S_{n-k}(u)$, where $2\le i \le d$, and $u$ is the center of $S_{n-k}$. If $T$ is a minimizing tree in $\T_n(d)$, then $T=\tilde{T}(i)$ for some $i$.
\end{theorem}

{\bf Conjecture 2.}
The maximizing tree in $\T_n(d)$ is $\tilde{T}(\lfloor d/ 2 \rfloor+1)$.

%\begin{con} Denote $P_{k+1}=v_{1}v_{2}\ldots v_{k+1}$, $\tilde{T}=P_{k+1}(v_{\lfloor (k+2)/ 2 \rfloor})\circ S_{n-k}(u)$, where $u$ is the center of $S_{n-k}$. Then for arbitrary tree $T$ with order $n$ and diameter $k$, $\rho(T^{2})\leq \rho(\tilde{T}^{2})$, the equality holds if and only if $T=\tilde{T}$.\end{con}

\end{document}